\documentclass [12pt,a4paper,reqno]{amsart}

\usepackage{amssymb}

 \input xy
\xyoption{all}

\usepackage{thmtools}
\declaretheoremstyle[notefont=\bfseries,notebraces={}{},%
    headpunct={},postheadspace=1em]{mystyle}
\declaretheorem[style=mystyle,numbered=no,name=Problem]{prob-hand}
\def\dispace{\setlength{\itemsep}{2pt}}

\newcommand{\onto}[1]{\;{\count255=0 \loop \relbar\joinrel
    \advance\count255 by1
    \ifnum\count255<#1 \repeat \twoheadrightarrow}\;}
\newcommand{\Onto}{\mathrel  - \joinrel \twoheadrightarrow}

\newcommand{\To}{\longrightarrow}

\def\iff{if and only if}

\def\lcm{\operatorname{lcm}}

\def\Dir{\Rightarrow}

\def\tl0{\widetilde 0}

\def\Q{\mathbb Q}

\newcommand{\etype}[1]{\renewcommand{\labelenumi}{(#1{enumi})}}
\def\eroman{\etype{\roman} \dispace}
\def\ealph{\etype{\alph}\dispace}

\def\pSkip{\vskip 1.5mm \noindent}
\newcommand{\ds}[1]{\ {#1} \ }
\newcommand{\dss}[1]{\quad {#1} \quad }

\def\semiring0{semiring$^{\dagger}$}

\def\lex{{\operatorname{lex}}}
\def\Fl{{\operatorname{Fl}}}
\def\Flmin{\Fl_{\min}}
\def\Fmin{F_{\min}}

\def\sm{\setminus}
\def\00{ \{ 0 \}}
\def\o00{\overline{\00}}
\def\onto{\twoheadrightarrow}

\def\|{\ds |}

\def\vrp{\varphi}

\def\Arch{{\operatorname{Arch}}}

\def\X1{X_1}
\def\Y1{Y_1}

\def\brV{\overline{V}}
\def\brR{\overline{R}}
\def\brS{\overline{S}}

\def\brC{\overline{C}}
\def\bra{\bar a}

\def\brx{\bar x}
\def\bry{\bar y}

\def\aV{{\operatorname{a} V}}
\def\aW{{\operatorname{a} W}}
\def\aU{{\operatorname{a} U}}

\def\conv{\operatorname{conv}}

\def\sat{{\operatorname{sat}}}

\def\N{\mathbb N}

\def\mfo{\mathfrak o}

\def\Aup{A^{\uparrow}}
\def\Zup{Z^{\uparrow}}

\def\Adw{A^{\downarrow}}

\def\Fdw{F^{\downarrow}}

\def\X1{X_1}
\def\Y1{Y_1}

\def\al{\alpha}

\def\Gm{\Gamma}

\newtheorem{thm}{Theorem} [section]
\newtheorem*{thm*}{Theorem}
\newtheorem{cor}[thm]{Corollary}
\newtheorem{lem}[thm]{Lemma}

\newtheorem{prop}[thm]{Proposition}

\newtheorem{rem}[thm]{Remark}
\newtheorem{rems}[thm]{Remarks}

\newtheorem*{claim*} {Claim}
\newtheorem*{theorem4.5'} {Theorem 4.5$'$}
\newtheorem{acknowledgment*}[thm] {Acknowledgment}

\newtheorem{dich}[thm]{Dichotomy}

\newtheorem{examp}[thm]{Example}

\newtheorem*{exampleA*}{Example A}
\newtheorem*{exampleB*}{Example B}
 \newtheorem{remark}[thm]{Remark}

 \newtheorem*{remark*}{Remark}
 \newtheorem{defn}[thm]{Definition}
\newtheorem{construction}[thm]{Construction}

\newtheorem{schol}[thm]{Scholium}

\newtheorem*{notation*} {Notation}
\newtheorem*{comment*} {Comment}

 \renewcommand{\sectionmark}[1]{}

\newcommand{\lm}{\lambda}
\newcommand{\Lm}{\Lambda}
\newcommand{\om}{\omega}

\newcommand{\On}{\operatorname{On}}

\begin{document}

\title[Archimedean classes in additive monoids]{Archimedean classes in additive monoids}
 \author[Z. Izhakian]{Zur Izhakian}
\address{
Department of Mathematics, Ariel University, 40700, Ariel, Israel.}
    \email{zzur@g.ariel.ac.il}

\author[M. Knebusch]{Manfred Knebusch}
\address{Department of Mathematics,
NWF-I Mathematik, Universit\"at Regensburg 93040 Regensburg,
Germany} \email{manfred.knebusch@mathematik.uni-regensburg.de}

\date{\today}







\begin{abstract} Summand absorbing submodules are common in modules over (additively) idempotent semirings, for example, in tropical algebra. A submodule $W$ of $V$ is  summand absorbing, if  $x + y
\in W$ implies  $x \in W, \; y \in W $ for any $x, y \in V$.  This paper proceeds  the    study of these submodules, and more generally of additive monoids,  with emphasis on their archimedean classes and quotient structures.   \end{abstract}

\maketitle

{ \small \tableofcontents}

\numberwithin{equation}{section}
\section*{Introduction}

The foundations for a module theory over semirings \cite{Cos,golan92} that follow the lines of classical module theory were laid in  \cite{Dec}. This theory has been developed further in \cite{SA,Gen,Aml}, and so far includes the notions of decomposition, generation,  amalgamation, and extension of summand absorbing modules. These notions help to bypass the lack of negation in modules over semirings and therefore play a significant role in the theory. The current paper proceeds the study of summand absorbing modules with emphasis on their archimedean classes and quotient structures, but in the more abstract category of additive monoids.

Aiming for a module theory over semirings, our initial motivation was to understand structures of modules in tropical algebra -- an algebra over a particular type of idempotent semirings. Nevertheless, there are many other important examples where such modules appear, e.g., additive semigroups, which can be viewed as modules over the semiring $\N_0$ of natural numbers, or sets of positive elements in an ordered ring or a semiring. Semirings appear extensively in other recent subjects of study (e.g. in discrete mathematics, logic,  and automata theory), providing a wide range of examples for modules over semirings.

The underpinning  attribute of summand absorbing modules is
 lack of zero sums:
An $R$-module $V$ over a semiring $R$ \textbf{lacks zero sums}
(abbreviated \textbf{LZS}),~ if
\begin{equation}\label{eq:LZS}
  \forall \, x, y \in V : \  x + y = 0 \Dir  x = y = 0. \tag{LZS}
\end{equation}
\noindent  Taking  submodules, direct sums, and direct products preserve LZS.
For example, the module $\operatorname{Fun} (S, V)$ of
 functions from a set $S$ to a
module~$V$ is LSZ \cite[Examples~1.6]{Dec}. More generally,   any  module over
an idempotent semiring is LZS \cite[Proposition~1.8]{Dec}, establishing a large assortment of
examples.

Relying on the property  LZS, a related type of submodules is obtained: A submodule $W$
of $V$  is \textbf{summand absorbing} (abbreviated \textbf{SA}) in
$V$, if
\begin{equation}\label{eq:SA}
\forall \, x, y \in V: \  x + y \in W \Dir x \in W, \; y \in W.
\tag{SA}
\end{equation} Such submodule $W$ is  called an \textbf{SA-submodule} of
$V$. In particular,  a left SA-ideal of a semiring $R$ is an SA-submodule of~$R$.

The two properties LZS and SA are linked: an $R$-module $V$ is LZS if and only if $\{ 0_V \} $ is an SA-submodule of $V$. This link enhances the interest in SA-submodules, although SA-submodules themselves are essential for any semiring $R$ and (left) $R$-modules~ $V$. Particular structures of SA-module have already appeared in the literature (e.g., see \cite{Tolliver}).

The merit of the SA-property to modules over a semiring was explored in our former papers, presenting applications and concrete examples \cite{Aml,Dec,SA,Gen}.  In particular,
\cite{Aml} reveals the symbiotic correspondence of the SA-property to amalgamations, a central technique in any module theory.
Our paper~ \cite{Dec} (jointly with L. Rowen) establishes the profound role of SA-modules as ``complements''  of submodules in a module that is   LZS.

Ignoring the scalar multiplication, every $R$-module is an additive monoid, for which  the definitions of properties  \eqref{eq:LZS} and \eqref{eq:SA} remain valid.
The paper takes this more general perspective  and focuses on an additive monoid $V$, which is \textbf{$D$-ordered}  by a submonoid $D \subset~ V$:
\begin{equation}\label{eq:Vqordering}
x \leq_D y \dss{\Leftrightarrow} x + d = y \text{ for some } d \in D. \tag{DO}
\end{equation}
This partial ordering allows for defining convexity,  archimedean property, and congruences.
When $D = V$, the additive  monoid $V$ with the ordering  $\leq_V$ is called  \textbf{upper bound} (or selective monoid \cite{golan92}).  Upper bound modules appear as  quotients of modules over semirings  by congruences, and are therefore associated with additive monoids, especially   with   SA-submonoids -- our main object of study.

Given an element $x > 0 $ in an additive monoid $V$, the \textbf{archimedean class} of $x$  is defined to be the set
\begin{equation}\label{eq:archClass}
 \Arch_V(x) := \{ y\in  V \ds | \exists m,n \in \N: x \leq_V ny, \ y \leq_V  mx \}.  \tag{AR}
\end{equation}
Namely, the so-called \textbf{axiom of Archimedes} holds, that assumes that any two lengths are comparable in the sense of \eqref{eq:archClass}, cf. Archimedes, \emph{On Spheres and Cylinders}, book I, Assumption 5,   the work
of Archimedes (Cambridge 1897) cited in \cite{Heath}. Ever since, it has been  an unsolved issue among historians, what Archimedes had in mind as the underlying structure for his famous study of lengths, volumes, and areas, e.g., see the list of books and papers in \cite[62f]{Struik}. The present paper suggests a much more abstract viewpoint to this idea of Archimedes, considering additive  monoids as underpinning structure.\footnote{
Some of the results in this paper can be phrased in terms of lattices. However, this paper focuses on monoid structures,
and leaves the lattice viewpoint to a future study.
}

The paper starts in \S\ref{sec:13} by assembling basic facts about archimedean classes to our framework, whose underlying structure is an \emph{arbitrary} additive monoid $V$ (or more specifically module). The relation $\leq_V$ is then only a quasiordering (reflexive, transitive, perhaps not antisymmetric) on an additive monoid $V$, utilized  to define an
\textbf{archimedean quasiordering}~ $\leq_\aV$. The quasiordering $\leq_\aV$  induces the equivalence relation $\equiv_\aV$  on $V$ that determines the  archimedean classes~ \eqref{eq:archClass}.

We show that, with respect to $\leq_\aV$, every archimedean class is convex in $V$ and closed for addition, and furthermore,  that the sum of two classes is contained in a unique class
(Theorem ~\ref{thm:1.3}). The latter property allows for defining an addition on the set $\Gm(V) $  of archimedean classes in $V$,  so that $\Gm(V)$ becomes a quotient  monoid together with the canonical projection
$\pi_\aV : V \twoheadrightarrow \Gm(V)$ that respects  quasiorderings  (Scholium \ref{schol:1.5}).
Then $\Gm(V) $ is an additively idempotent monoid, i.e., $\xi + \xi = \xi$ for every $\xi \in \Gm(V) $.
This setup  results in a homomorphism theorem of additive monoids and their quotient monoids  (Theorem \ref{thm:1.9}).

Employing the quasiordering $\leq_V$, each additive subset $S$ of  $V$, i.e., $S+S \subset S$, determines an SA-submonoid $D = \brC(S)$  (cf. \S\ref{sec:13}).
 Such $S$ is called an  \textbf{attractor} of $D$ (Definition ~\ref{def:2.3}).
 When $S = \N x$, the SA-submonoid $ \brC(\N x )$ is precisely the archimedean class $\Arch_V(x)$.
The submonoid  $\brC(S)$ is the union of all archimedean classes it contains. In particular, $\brC(A) = A$ for any archimedean class~ $A$, and the sum of two archimedean classes is contained in a unique archimedean class.

 We further have a natural notion of \textbf{saturation} of  any additive subset $S$ in $V$ (Definition~ \ref{def:2.4}).
 Theorem ~\ref{thm:2.9}  links a minimal SA-submonoid $W(A)$ containing a given  archimedean class $A$ in ~ $V$ to a saturated additive set and characterizes this set. Theorem~ \ref{thm:2.10} translates the relation between classes $A, B \in \Gm(V)$ to the corresponding SA-submonoids $W(A)$,~ $W(B)$.

A \textbf{flock} $\Fl(A)$ in an additive monoid $V$ is a set $S \subset \Gm(V)$ of  archimedean classes, closed for addition,  which are compatible to the given class $A$.  Minimality of classes  is  discussed in detail at the end of  \S\ref{sec:2}.
Theorem \ref{thm:3.2} shows that the SA-submonoid $W(A)$ is a disjoint union of the minimal classes in $\Fl(A)$. Construction \ref{cons:3.8} defines a sequence of subsets of a maximal flock $F$, which is unique.  We conclude that, every $A \in \Gm(V) $ is a member of a unique maximal flock.
Then~ \S\ref{sec:3} analyzes the decomposition of a maximal flock $F$ to smaller flocks, and links these  to a height function on $F$, with  values in ordinal numbers $ < \om + \om$.

Our study of a height function in \S\ref{sec:3} establishes a hierarchy of archimedean classes in maximal flocks, based on the specific  ``scale'' $(\N,+)$.
Our motivation is to employ more general scales $(\Lm, +)$ of totaly ordered semigroups, which would open the area of archimedean phenomena to other techniques, in particular mathematical logic and
model theory, and their applications.  In the present paper, we include only a primordial case of this idea, aiming to give the reader a view  to this vision.

\section{Preliminaries}\label{sec:13}
To make the paper reasonably self-contained, we recall some basic definitions and relevant results,  mostly from  \cite[\S13]{Aml}. Recall that any $R$-module, is a an additive monoid, so that properties  of an additive monoid hold also for any $R$-module. In what follows, unless otherwise specified, $V = (V,+)$ stands for an additive monoid.

The $V$-quasiordering $\leq_V$ on an additive  monoid $V$, cf. \eqref{eq:Vqordering} with $D = V$, determines the congruence relation  $x \equiv_V y$ on $V$, given by
$$ x \equiv_V y \dss\Leftrightarrow x \leq_V y \text{ and } y \leq_V x.$$
Dividing out the monoid $V$ by the congruence  relation  $\equiv_V $ gives the monoid $\brV = (\brV,+)$, where $\equiv_V $  turns the quasiordering  $\leq_V$ on $V$ to a partial ordering $\leq_{\brV}$ on $\brV = V / \equiv_V$.
We denote  by $\brx$ the congruence class of an element  $x \in V$, and have for $x,y \in V$ the explicit description
\begin{equation}\label{eq:13.1}
  \brx = \bry \dss\Leftrightarrow \exists z,w \in V: x+z = y,\ y+w = x.
\end{equation}
The monoid $\brV$ is called the \textbf{upper bound  monoid associated to $V$}.

\begin{remark} Given a semiring $R$, the equivalence relation $ \equiv_R$ is a semiring congruence. Indeed,
if $a_1 \equiv_R a_2 $ and $b_1 \equiv_R b_2 $, then
$a_1 +d_1 = a_2 $, $a_2 + d_2 = a_1$,  $b_1 +e_1 = b_2$, $b_2 +e_2 = b_1$.
Then, for the multiplication,  $a_2 b_2 = (a_1 +d_1) (b_1 +e_1) = a_1 b_1 + (e_1 a_1  + d_1 b_1  + d_1 e_1)$, implying that $a_1 b_1 \leq_R a_2 b_2$. By symmetry $a_2 b_2 \leq_R a_1 b_1$, and thus
$a_1 b_1 \equiv_R a_2 b_2$.
For the addition,  $a_2 + b_2 = (a_1 +d_1) + (b_1 +e_1) = a_1 + b_1 + ( d_1 + e_1)$, implying that $a_1 +  b_1 \leq_R a_2 + b_2$. By symmetry $a_2 + b_2 \leq_R a_1 + b_1$, and thus
$a_1 + b_1 \equiv_R a_2 + b_2$.
\end{remark}

If $V$ is a module   over a semiring $R$, then  $\brV$ is a module over the upper bound semiring $\brR = R / \equiv_R$, with scalar mutiplication
$  \bra \bar v = \overline{av}
$ for $a \in R $, $v \in V$. Hence
$\sum_{i=1}^{r}  \bra_i {\bar v}_{i} = \overline{ \sum_{i =1}^r a_i v_i } $ for $a_i \in R $, $v_i \in V$.
  Then  $\brV$ is called the upper bound  $\brR$-\textbf{module associated to $V$}.
Most of our work is done in the category of $R$-modules, where   $R$ is an upper bound semiring. Monoids  can be subsumed here  by taking
$R = \brR =\N_0$.

\begin{defn}
A subset $S \subset V $ is \textbf{convex}, if $s \in S$  for any $s_1, s_2 \in S$, and $s\in V$ with $s_1 \leq_V s \leq_V s_2$. The \textbf{convex hull} of a subset $T \subset V$ is denoted by $\conv(T)$.
\end{defn}

Denote the convex hull of $\N_0 \cdot 1_R$ in $R $ with respect to $\leq_R$ by
$$\mfo_R = \conv(\N_0\cdot 1_R).$$
\begin{remark} The convex hull
$\mfo_R = \conv(\N_0\cdot 1_R)$ is a semiring. Indeed, clearly $0 \cdot 1_R = 0_R \in \mfo_R $ and
$1 \cdot 1_R = 1_R \in \mfo_R $. If $a, b \in \mfo_R $, then $m_1 \cdot 1_R \leq_R a \leq_R m _2 \cdot 1_R $ and
$n_1 \cdot 1_R \leq _R b \leq _R  n_2 \cdot 1_R $ for some $m_1,m_1,n_1, n_2 \in \N_0$. Then
$(m_1 + n_1) \cdot 1_R   \leq_R  a+ b \leq_R  (m_2+n_2) \cdot 1_R$ and
$(m_1  n_1) \cdot 1_R   \leq_R  a b \leq_R  (m_2 n_2) \cdot 1_R$. Hence $a+ b \in \mfo_R $ and
$ab \in \mfo_R $.
\end{remark}

Often is it better to work with a module $(V,+)$ over a semiring $R$ instead of an additive monoid. Modules over semirings restrict to additive monoids by taking  $R= \N_0$, and have a richer structure, where the multiplicative  monoid  $(R, \cdot \, )$ servers as a set of operators on the $R$-module $ (V,+)$. Then, $\mfo_R$ is realized as a semiring by applying this view to  $V=R$ with $R$ operating on $R$ by multiplication.

The following useful property holds for  any monoid $(V,+)$, and thus for any module $V$
\begin{prop}[{\cite[Lemma 6.6]{Dec}}]\label{prop:13.1}
  A submodule $S$ of a module $V$ is SA in $V$ \iff \ $S$ is a union of equivalence classes in $\brV = V/ \equiv_V$ and $\brS= S/ \equiv_V $ is SA in $\brV$.
\end{prop}

Assuming that  $V$ is an upper bound monoid, for every $x \in V$ we define the SA-monoid $$C(x) = \{ v \in V \ds | v+x =x\} $$ of $V$  \cite[\S12]{Aml}. Concerning  the upper bound monoid $\brV = V / \equiv_V$ associated to \emph{any} additive monoid $V$,   by the use of Proposition \ref{prop:13.1}, instead of arguing in $\brV$ and then passing to~ $V$, we work directly in $V$ with the quasiordering $\leq_V$.  Given $x \in V$ we define
\begin{equation*}\label{eq:13.3}
  \brC(x) := \{ u \in V \ds | u +x \leq_V x\},
\end{equation*}
 which means  (since $x \leq_V u +x$)
\begin{equation}\label{eq:13.4}
  \brC(x) = \{ u \in V \ds | u +x \equiv_V x\}.
\end{equation}
Thus $\brC(x) = C(x)$, when $V$ happens to be upper bound.

Assuming that $V$ is an $R$-module, where $R$ is any semiring, then $\brV$ is an $\brR$-module for $\brR = R/\equiv_R$ and $\brC(x)$ is an $\mfo_R$-submodule  of $V$  \cite[Corollary 12.10]{Aml}.

\pSkip

We write $nx$, $n \in \N$,  for the sum $x+ \cdots  +x$ with $x$ repeated  $n$ times.

\begin{prop}[{\cite[Propositions 13.2 and 13.3]{Aml}}]\label{prop:13.2}
$ $
\begin{enumerate}\eroman
  \item $\brC(x)$ is an SA-submonoid  of $V$.

  \item If $x \leq_V x'$ (resp.  $x \equiv_V x'$), then $\brC(x) \subset \brC(x')$ (resp.  $\brC(x) = \brC(x')$).
\end{enumerate}
\end{prop}
Therefore
\begin{equation}\label{eq:13.5}
   \brC( n x)  \subset \brC(nx + x) = \brC((n+1)x).
\end{equation}
Clearly,   the subset
$$ \brC_\om(x) := \bigcup_{n \in \N} \brC(nx)$$
 is again an SA-submonoid of $V$, which can be  defined as
$ C_\om(x) := \bigcup_{n \in \N} C(nx)$  in the case that $V$ is upper bound.
\begin{defn}[{\cite[Definition 13.5]{Aml}}]\label{def:13.5} The \textbf{archimedean class} $\Arch_V(x)$ of an element $x \in V$ is the set of all $y \in V$ such that $x \leq_V n y$ and $y \leq_V m x $ for some $n,m \in \N$.
\end{defn}

\begin{prop}[{\cite[Lemma 13.6, Proposition 13.7]{Aml}}]\label{prop:13.7}
$ $

\begin{enumerate}\eroman
  \item Assume that $x \leq_V my$ for some $m \in \N$.
Then $\brC(x) \subset \brC(my)$ and $\brC_\om(x) \subset \brC_\om(y)$.

  \item  If $\Arch_V(x) = \Arch_V(y)$, then $\brC_\om(x) = \brC_\om(y)$.

\end{enumerate}
 \end{prop}
Given a subset $S$ of $V$ such that  $S +S \subset S$, the subset
\begin{equation}\label{eq:13.6}
   \brC(S) := \bigcup_{s \in S} \brC(s)
\end{equation}
of $V$ is an SA-submonoid. In particular, $\brC(s+t)$ is an SA-submonoid of $V$, since
$$ x + y \in  \brC(s) +  \brC(t) \subset  \brC(s + t), \qquad x \in \brC(s), y \in \brC(t).$$
  Again $\brC(S)$ is an $\mfo_R$-submodule of $V$, when  $V$ is an $R$-module. Also
\begin{equation}\label{eq:13.7}
   \brC_\om(x) =  \brC( \N x) \qquad \text{for any $x \in V$.}
\end{equation}


Given two subsets $S$ and $T$ of $V$, closed under addition, suppose that for any $t \in T$ there exists $s \in S$ with $t \leq s$. Then $\brC(T) \subset \brC(S)$. \emph{It follows that $\brC(T) = \brC(S)$, if $S$ and $T$ are cofinal under $\leq_V$.}

\section{The monoid of archimedean classes}\label{sec:1}

In this section, we focus on archimedean classes of an additive monoid $V$ in a more comprehensive and systematic way than in \cite[\S 12 and \S13]{Aml}.

Starting with the quasiordering $\leq_V$ (often written simply as $\leq$),  cf. \eqref{eq:Vqordering}, we introduce the following coarser quasiordering $\leq_\aV$ on~ $V$, named \textbf{archimedean quasiordering}:
\begin{equation}\label{eq:1.1}
  \begin{array}{lll}
    x \leq_\aV x' & \Leftrightarrow & \exists n \in \N : x \leq_V n x' \\ [1mm]
    &  \Leftrightarrow & \exists n \in \N , z \in V: \  x + z =  n x'. \\
  \end{array}
  \end{equation}
  Consequently, we have the associated equivalence relation $\equiv_\aV$, given by
\begin{equation}\label{eq:1.2}
  \begin{array}{lll}
    x \equiv_\aV x' & \Leftrightarrow & x \leq_\aV x', \ x' \leq_\aV x. \\
  \end{array}
  \end{equation}
The equivalence classes of  $\equiv_\aV$ are the archimedean classes, as defined in
Definition \ref{def:13.5},
 where  now we denote the archimedean class $\Arch_V(x)$ of $x$ by $[x]_\aV$, for short.

Observe that the quasiordering $\leq_\aV$ is compatible with the addition of $V$:
\begin{equation}\label{eq:1.3}
  \begin{array}{lll}
    x \leq_\aV x', \   u \leq_\aV u' & \Rightarrow & x + u \leq_\aV x' + u'. \\
  \end{array}
\end{equation}
Indeed, if $x \leq_V n x'$ and $u \leq_V m u'$, and say $n \leq_V m,$ then $x+ u \leq_V m (x' + u')$.
By \eqref{eq:1.3} it follows that
\begin{equation*}\label{eq:1.4}
  \begin{array}{lll}
    x \equiv_\aV x', \ u \equiv_\aV u' & \Rightarrow & x + u \equiv_\aV x' + u'. \\
  \end{array}
\end{equation*}
Furthermore,  from \eqref{eq:1.1} and \eqref{eq:1.2} we conclude that
\begin{equation}\label{eq:1.5}
  \forall n \in \N, x \in V: \quad nx \equiv_\aV x.
\end{equation}
The zero class $[0]_\aV$ is the set of all $y \in V$ with $0 \leq_V my$ and $y \leq_V n 0$ for some $m,n \in \N$. The former condition holds for all $y \in V$, and the latter means that $y \leq_V 0$.
Thus
\begin{equation*}\label{eq:1.6}
  [0]_\aV = \{ y \in V \ds | \exists z \in V : \ y + z = 0\}.
\end{equation*}
In other terms, cf. \eqref{eq:13.4},
\begin{equation}\label{eq:1.7}
  [0]_\aV = \brC(0).
\end{equation}

In the study of archimedean classes it will be of help to use a type of scalar multiplication of convex subsets of $V$, which are closed under addition, by suitable elements of $\Q$. For any (nonempty) subset $A$ of $V$ we define
\begin{equation*}\label{eq:1.8}
  \frac 1n A = \frac 1 n \cdot_V A := \{ x\in V \ds | nx \in A \}.
\end{equation*}
On the other hand, we use the notation
\begin{equation*}\label{eq:1.9}
  m  A := \underbrace{A + \cdots + A}_m
\end{equation*} for any $m \in \N$.
It is obvious that for any $n \in \N$
\begin{equation}\label{eq:1.10}
  A \subset \frac 1n A
\end{equation}
and for any $n_1, n_2 \in \N$
\begin{equation}\label{eq:1.11}
  \frac{1}{n_1}\bigg(\frac{1}{n_2} A\bigg) = \frac{1}{n_1 n_2} A.
\end{equation}
\begin{lem}\label{lem:1.1}
If $A$ is convex in $V$, then $\frac 1n A$ is also  convex in $V$.
\end{lem}
\begin{proof}
  Let $x, x' \in \frac 1n A$, $y \in A$, where $x \leq_V y \leq_V x'$. Then $nx \leq_V n y \leq_V nx'$ and $nx, nx' \in A$. Since~ $A$ is convex, it follows that $n y \in A$, whence $ y \in \frac 1n A$.
\end{proof}

\begin{lem}\label{lem:1.2}
Assume that $A$ is closed under addition. Then for any $n_1, n_2 \in \N$
$$ \frac{1}{n_1}A + \frac{1}{n_2}A \subset \frac{1}{n} A,$$
where  $n = \lcm(n_1, n_2)$.
\end{lem}
\begin{proof}
A straightforward consequence of \eqref{eq:1.10} and \eqref{eq:1.11}.
\end{proof}
\begin{thm}\label{thm:1.3} $ $
\begin{enumerate} \ealph
  \item Every archimedean class $A$ of an additive monoid $V$ is a convex set in $V$,  closed under addition, and satisfies $A = \frac 1n A$ for any $n \in \N$.
  \item For any two  archimedean classes $A$ and $B$ the sum
  $$ A + B := \{ x + y \ds | x \in A, y \in B\} $$ is contained in a unique  archimedean class $C$.
  \item The class $C$ contains the union
  \begin{equation*}\label{eq:1.12}
    C'(A,B) := \bigcup_{n \in \N} \frac 1n \conv(A+B),
  \end{equation*}
  where $\conv(A+B)$ denotes the convex hull of $A+B$. This union  is convex and closed under addition.
\end{enumerate}

\end{thm}
\begin{proof}
(a):
  Given two archimedean classes $A = [x]_\aV$ and $B = [y]_\aV$, it is clear by \eqref{eq:1.3} that $A  + B= [x+y]_\aV$.
  In the case $A = B$, we further have
  $[x]_\aV = [x+ x]_\aV$ (cf. \eqref{eq:1.5}), and this class is $A$ itself. Thus  $A$ is closed under addition, and  convex. Indeed, if $x_1, x_2 \in A$ and $x_1 \leq_V  y \leq_V x_2$, then
$x_1 \leq_V  y \leq_V x_2 \leq _\aV x_1$, and thus
$x_1 \leq_\aV  y \leq_\aV x_1$, whence $[x_1]_\aV = [y]_\aV$.
\pSkip
(b):
Applying (a) to the classes $A = [x]_\aV$ and $B = [y]_\aV$, we see that $A+B$ is closed under addition, and the unique archimedean class
$C = [x+y]_\aV $ containing $A+B$ is convex, whence $\conv(A+B) \subset C$.
\pSkip
(c):
By use of  \eqref{eq:1.5},  we conclude that $C'(A+B) \subset C$, and by Lemma \ref{lem:1.2} it is now evident that $C'(A,B)$ is convex and closed under addition.
\end{proof}

\begin{defn}\label{def:1.4}
Let $\Gm(V)$ denote the set of all archimedean classes of an additive monoid~ $V$. Given $A,B \in \Gm(V)$, we define the sum $A +_\Gm B$ in $\Gm(V)$ to be the unique class $C \in \Gm(V)$ for which  $ A+B \subset C$, cf. Theorem \ref{thm:1.3}.(b).
\end{defn}
\noindent The set $\Gm(V)$ is an idempotent additive monoid, since clearly $A+ A = A$.
\pSkip

We aim to realize  the assignment $V \mapsto \Gm(V)$ as a functor in the category of additive monoids.

\begin{schol}\label{schol:1.5}  We consider the archimedean quasiordering $\leq_\aV$ on $V$, and the associated equivalence relation $\equiv_\aV$, given by \eqref{eq:1.1}  and  \eqref{eq:1.2}.
Both are compatible with addition. Then $\Gm(V)$ is the quotient monoid
$V / \equiv_\aV$ with the canonical projection
\begin{equation*}\label{eq:1.13}
  \pi_\aV: V \Onto \Gm(V), \qquad x \mapsto [x]_\aV
.\end{equation*}
The monoid
 $\Gm(V)$ is equipped with  the partial ordering
\begin{equation*}\label{eq:1.14}
[x]_\aV \leq_\Gm [y]_\aV \dss \Leftrightarrow x \leq_\aV y,
\end{equation*}
and the  addition
\begin{equation}\label{eq:1.15}
[x]_\aV +_\Gm [y]_\aV = [x+y]_\aV ,
\end{equation}
both  are well defined. 
\end{schol}
We state some basic facts about idempotent monoids.
\begin{lem}\label{lem:1.6} Assume that $V$ is an idempotent monoid, i.e., $x+ x = x$ for all  $x \in V$.
\begin{enumerate}\ealph \item Then  $V$ is upper bound and, for any $x, y \in V$,
$$y \leq_V x  \dss \Leftrightarrow x +y= x.$$

  \item For any $x, x' \in V$,
$$x\leq_\aV x'  \dss \Leftrightarrow x \leq_V x'.$$

\end{enumerate}
\end{lem}
\begin{proof} (a): The relation
$y \leq_V x$ means that there exists $z \in V$ such that  $y + z =x$. Thus  $y +z =x \Rightarrow x+ y +z = x \Rightarrow x+y =
(x+ y +z + y) = x$. This proves both claims.
\pSkip
(b): The relation
  $x\leq_\aV x'$ means that $x \leq_V n x'$ for some $n \in \N$. In the present case
  $nx' =x'$.
\end{proof}
Lemma \ref{lem:1.6} gives us a complement to the preceding scholium.
\begin{schol}\label{shcol:1.7} The following properties are equivalent for $V$.
\begin{enumerate} \ealph
  \item $[x]_\aV = [x]_V$ for all $x \in V$,
  \item $\pi_\aV: V \to \Gm(V)$ is a monoid isomorphism.
\end{enumerate} Properties
(a) and (b) hold, if $V$ is idempotent.
\end{schol}

\begin{lem}\label{lem:1.8}  Let $\vrp: V \to W$ be a homomorphism  from $(V,+)$ to an additive monoid ~$(W,+)$.
\begin{enumerate}\ealph
  \item If $x,x' \in V$ and $ x\leq_\aV x'$, then $\vrp(x) \leq_\aW \vrp(x')$.
  \item  If $ x \equiv_\aV x'$, then $\vrp(x) \equiv_\aW \vrp(x')$.
  \item For any archimedean class $A \in \Gm(V)$ there exists a unique class $B \in \Gm(W)$ such that $\vrp(A) \subset B$.
\end{enumerate}
\end{lem}
\begin{proof} The relation
  $x \leq_\aV x'$ means that there exist $n \in \N$ and $y \in V$ with $x+y = nx'$. This implies $\vrp(x) + \vrp(y)= n \vrp(x') $. Thus (a) holds true,  and then (b) and (c) follow.
\end{proof}

From these  considerations we infer the following in a straightforward way.
\begin{thm}\label{thm:1.9} $ $
\begin{enumerate} \ealph
  \item Given a homomorphism $\vrp: V \to W$ of additive monoids, there exits a well-defined monoid homomorphism
  $$ \Gm(\vrp): \Gm(V) \To \Gm(W), \qquad [x]_\aV \mapsto [\vrp(x)]_\aW,$$
  which maps each class $A \in \Gm(V)$ to the unique class $B \in \Gm(W)$  containing  $\vrp(A)$.
The diagram \begin{equation*}
    \xymatrix{
    V    \ar@{->}[rr]^{\pi_\aV} \ar@{->}[d]^{\vrp} & & \Gm(V)  \ar@{->}[d]^{\Gm(\vrp)} \\
   W  \ar@{->}[rr]^{\pi_\aW}   &  &  \Gm(V)
   }
  \end{equation*}
commutes.

\item If the monoid $V$ is idempotent, then $\pi_\aV$ is an isomorphism.
\end{enumerate}
\end{thm}

It is easily seen that the functor $\Gm$ respects arbitrary direct sums.  We explore this property  in the case of two summands.
\begin{examp}\label{exmp:1.10}
Assume that $V_1$ and  $V_2$ are additive monoids and $V = V_1 \oplus V_2 = V_1 \times V_2$. If $x_1, x_1' \in V_1$,
$x_2, x_2' \in V_2$, then
$$ x_1 + x_2 \leq_V x_1' +x_2' \dss \Leftrightarrow x_1 \leq_{V_1} x_1', \ x_2 \leq_{V_2} x_2'. $$
  Thus,  by  \eqref{eq:1.1} and \eqref{eq:1.2}  we conclude that,
\begin{equation*}\label{eq:1.16}
    x_1 + x_2 \leq_\aV x_1' +x_2' \dss \Leftrightarrow x_1 \leq_{\aV_1} x_1', \ x_2\leq_{\aV_2} x_2',
  \end{equation*}
for any $x_1, x_1' \in V_1$,
  $x_2, x_2' \in V_2$, and (consequently)
  \begin{equation*}\label{eq:1.16}
    x_1 + x_2 \equiv_\aV x_1' +x_2' \dss \Leftrightarrow x_1 \equiv_{\aV_1} x_1',\  x_2\equiv_{\aV_2} x_2'.
  \end{equation*}
  Thus, as  partially ordered semigroups, we have
  $ \Gm(V) = \Gm(V_1) \times \Gm(V_2)$.
  Observing that, in this notation,
  $ [0]_\aV = ([0]_{\aV_1}, [0]_{\aV_2})$,
we obtain that
\begin{equation*}\label{eq:1.18}
   \Gm(V) = \Gm(V_1) \oplus \Gm(V_2)
\end{equation*}
in the category of ordered monoids.
\end{examp}
\begin{examp}\label{exmp:1.11}
Given two  upper bound additive monoids $V_1$ and $V_2$, we equip the product of sets $V = V_1 \times V_2$ by  the lexicographic ordering
\begin{equation*}\label{eq:1.18}
   V = V_1 \times_\lex V_2,
\end{equation*}
where a tuple $(x_1, x_2)$ precedes a tuple $(x_1', x_2')$, if either $x_1 \leq_{V_1} x_1'$, or  $x_1 = x_1'$ and
 $x_2 \leq_{V_2} x_2'$. Assuming  that $V_1$ is cancellative,  the ordering $\leq_V$ on $V$ is compatible with addition, and thus ~ $V$ is an upper bound additive monoid.

 Given $x_1, x_1' \in V_1$, $x_2, x_2' \in V_2$,  we have
 $(x_1, x_2) \equiv_\aV (x_1', x_2')$ \iff \ either $x_1 \neq 0$ and $x_1 \equiv_{\aV_1} x_1'$, or $x_1 = x_1' =0$ and $x_2 \equiv_{\aV_2} x_2$.
  Thus $\Gm(V)$ is the union of the two convex subsets $A\cong \Gm(V_1) \sm \00$ and  $B \cong \Gm(V_2)$, where $B < A$.
\end{examp}

\begin{examp}\label{exmp:1.12}
Assume that $V_1$ is an idempotent monoid, $V_2$ an additive monoid, and take  the lexicographic product
$V = V_1 \times_\lex V_2$. Then
$ \Gm(V) = V_1 \times \Gm(V_2)$
is the disjoint union of ~$\00$ and of copies $\{ x\} \times \Gm(V_2)$ of $\Gm(V_2)$, which are incomparable.
\end{examp}

\section{The SA-submonoids $\brC_\om(x)$ and their attractors}\label{sec:2}

 SA-submonoids $\brC_\om(x)$ are the main objects of study in this section. We begin with basic characterizations of certain subsets.
\begin{defn}\label{def:2.1} $ $
\begin{enumerate} \ealph
\item
We call a subset $D \subset V$ an \textbf{entourage} in $V$, if there exists some $x \in V$, for which $\brC_\om(x) = D$. Recall from  \eqref{eq:13.6} that then $D$ is an SA-submonoid of $V$.

\item We call a set $S \subset V$, which is closed under addition, an \textbf{attractor} of $D$, if \ $\brC(S) =~ D$, and then also say  that $S$ \textbf{attracts} $D$. We denote the set  of attractors of $D$ by~ $T(D)$.

\end{enumerate}
\end{defn}
\noindent In this terminology $\brC_\om(x) = D$ means that $\N x$ is an attractor of $D$, cf.
\eqref{eq:13.7}.
\pSkip

If $S$ is an attractor of $D$, and $S'$ is a subset of $V$ such that  $S' + S' \subset S'$, which is cofinal to~ $S$, then clearly $S'$ is also an attractor of $D$. Thus all cofinal sets $S$ to $\N x$ with  $S+S \subset S$ are attractors of $D$. An archimedean class $[x]_\aV$ is such a maximal  subset of ~$V$, as is clear from the definition of $\equiv_\aV$,
$$ x \equiv_\aV x' \dss \Leftrightarrow \exists m,n \text{ with } x \leq_V nx', \ x' \leq_V mx.$$
Hence, when looking for attractors of an entourage, we may focus on archimedean classes.

\begin{dich}\label{dich:2.2}
  Given a subset $D$ and an archimedean class $A$ that attracts $D$, we face two cases.

  \begin{description}
    \item[Case I] There exists $z \in A \cap D$. Then such $z$ is unique up to $\equiv_\aV$. Indeed, if also $z' \in A \cap D$, then $z+ z \leq_V z$, $z+ z' \leq_V z'$, $z'+ z' \leq_V z'$, and thus
        $ z \equiv_V z+z \equiv_V z+ z' \equiv_V z' + z' \equiv_V z'$. In short $A = [z]_\aV$, and so $A$ is a singleton up to
        $\equiv_\aV$.
        We then call $A = [z]_\aV$ (or abusively the element $z$) the \textbf{center} of $D$. \pSkip
    \item[Case II] $A \cap D = \emptyset$. Then, for any $x \in A$,
    $$ \brC(x) \ds \subsetneqq  \brC(2x) \ds \subsetneqq  \brC(3x) \ds \subsetneqq \cdots, $$
    and $D = \brC_\om(x)$.
  \end{description}
\end{dich}

We next  explore  subsets $S$ of $V$ having the property $S+S \subset S$, which can serve as attractors of a given entourage $D$ in $V$.
\begin{defn}\label{def:2.3} $ $
\begin{enumerate} \ealph
\item A \textbf{saturated additive} subset $S$ of $V$ is a set $S \subset V$ satisfying  $S+S \subset S$, and closed under the equivalence $\equiv_\aV$. In other terms, $S$ is the preimage of a set $\brS \subset \Gm(V) = V/ \equiv_\aV$ such that  $\brS + \brS \subset \brS$.

    \item A subset $S$ is called a \textbf{principal additive subset} of $V$, if $S$ is saturated additive and there exists  $x \in S$ for which
    $\brC(S) = \brC_\om(x)$.

\item
     If $S$ is a principal additive subset, and $x \in S$, then either
     $\brC(\N x) = \brC([x]_\aV) = \brC(S)$, or $\brC(\N x) \subsetneqq \brC(S)$.
In the first case we say that the element $x$ or the archimedean class     $A = [x]_\aV$ is \textbf{essential} in $S$. In the second case we call $x$ or $A$ \textbf{accessory} (or \textbf{inessential}) in $S$.
\end{enumerate}
\end{defn}

\begin{defn}\label{def:2.4} Given an additive subset $X$ of $V $ (i.e., $X+X \subset X$), we call the union $X_\sat$ of all classes $[z]_\aV$, $z \in X$, the \textbf{saturation} of $V$.\footnote{The saturation of $V$ can be realized as a closure
operator on the additive subsets of $V$.}
\end{defn}

 The set $X_\sat$  is again additive, since
$$[x]_\aV + [y]_\aV = [x+y]_\aV$$
for $x,y \in X$, and thus is the smallest saturated additive set in $V$ that contains $X$. \pSkip

We state some facts, which are clear by the previous definitions and considerations.
\begin{rems}\label{rem:2.5}
Let $\vrp: V \to W$ be any homomorphism of additive monoids.

\begin{enumerate} \ealph
\item For any two additive subsets $X,Y$ of $V$ we have
$$ X+Y \subset X_\sat + Y_\sat \subset (X+Y)_\sat.$$

\item If $A_1$ and $A_2$ are archimedean classes in $V$, then (cf. \eqref{eq:1.15})
$$ A_1 +_\Gm A_2 = (A_1+ A_2)_\sat.$$

\item Let  $A = [x]_\aV \in \Gm(V)$.  Then
$$ B := \vrp(A)_\sat = [\vrp(x)]_\aW$$
is the unique archimedean class of $W$ that contains $\vrp(A)$, and $B = \Gm(\vrp(A))$.

\item If $D$ is an  entourage in $V$ and $A \in \Gm(V)$ is an attractor of $D$, where  $D = \brC(A)$, then $E = \vrp(D)_\sat$ is an entourage in $W$ and $B = \vrp(A)_\sat$ is an attractor of $E$. Recall that ~$D$ and $E$ are SA-submonoids of $V$ and $W$ respectively.
\end{enumerate}
\end{rems}
Observe that, if $A_1$ and $A_2$ are archimedean classes attracting a given entourage $D$, then either
$\brC(A_1 + A_2) = D$ or  $\brC(A_1 + A_2) \supsetneqq D$. The first case reflects a kind of ``sympathy'' between the classes $A_1$ and $A_2$, which is missing in the second case. This observation prompts an equivalence relation on the set
$\Gm(V)$ of archimedean classes as follows.

\begin{defn}\label{def:2.6} $ $
\begin{enumerate} \ealph
\item We call two archimedean classes $A_1, A_2 \in \Gm(V)$ \textbf{relatives}, or say that $A_1$ and $A_2$ are \textbf{compatible}, if
    \begin{equation*}\label{eq:2.1}
      \brC(A_1) = \brC(A_2) = \brC(A_1 + A_2).
    \end{equation*}
    Since the additive set $A_1 + A_2$ is contained in a unique archimedean class $A_1 +_\Gm A_2$, cf.  Definition \ref{def:1.4}, this means that
      \begin{equation*}\label{eq:2.2}
      \brC(A_1) = \brC(A_2) = \brC(A_1 + _\Gm  A_2).
    \end{equation*}

    \item A \textbf{flock} in $V$ is a (nonempty) set $F$ of compatible archimedean classes in $V$, which is closed under addition in $\Gm(V)$. This means that, if $A_1, A_2 \in F$, then
        $$\brC(A_1+A_2) = \brC(A_1)+\brC(A_2)$$ and $A_1 +_\Gm A_2 = (A_1 + A_2)_\sat \in F.$
        Thus, there exists a unique entourage $D \subset V$ with $D = \brC(A)$ for every $A \in F$. We say that $F$ is a \textbf{flock of attractors} of $D$.

\end{enumerate}
\end{defn}
It is  obvious that every $A \in \Gm(V)$ is a member of a unique \textbf{maximal flock}, i.e.,  the ~set
\begin{equation}\label{eq:2.3}
  \Fl(A) = \{ A' \in \Gm(V) \ds | \brC(A) = \brC(A') = \brC(A+A')    \}.
\end{equation}

\begin{prop}\label{prop:2.7}
  Let $W$ be an SA-submonoid of $V$, and let $A$ be an archimedean class in~ $V$.
  \begin{enumerate} \ealph
    \item If $A \cap W \neq \emptyset$, then $A \subset W$. If $D$ is the entourage of $V$ with attractor $A$, then
    $D \subset W$, and so is the entourage of $W$ attracted by $A$.
    \item  $\Gm(W) \subset \Gm(V)$, and
    $$ \Gm(V) \sm \Gm(W) = \{ B \in \Gm(V) \ds | B \cap W = \emptyset\}. $$
  \end{enumerate}
\end{prop}
\begin{proof} (a): Choose some $x \in A \cap W$. If $x' \in A$, then there exist $n \in \N$ and $y \in V$ for which
$n x = x' +y$. Since $nx \in W$ and $W$ is SA in $V$, this implies that $x' \in W$ (and $y \in W$), and proves that $A \subset W$. If $d$ is contained in the entourage $D$ attracting  $A$, then there exists   $n \in \N$ for which
$$ d + nx \equiv_V nx \in A \subset W,$$
   implying  that $d \in W$, an thus $D \subset W$. \pSkip
  (b): Now obvious.
\end{proof}

It follows by Proposition \ref{prop:2.7}.(a) that for a given archimedean class $A$ in $V$ there exists a unique minimal summand absorbing submonoid $W$ of $V$ that  contains the set $A$, and this submonoid~ $W$ is saturated in $V$ (i.e., a union of archimedean classes).
We search for a description of~ $W$.
To warm up, we state an easy observation about saturated submonoids related to $A$, not necessarily summand absorbing.
\begin{rem}\label{rem:2.8} Let $\o00 = [0]_\aV = \brC(0)$, i.e., the smallest summand absorbing submonoid of ~$V$, cf.~ \eqref{eq:1.7}.  Then $A \cup \o00$ is the smallest saturated submonoid of $V$ containing $A$ -- the submonoid generated by $A$. If $A \cap \brC(A) = \emptyset$, then
    \begin{equation}\label{eq:2.4}
    \begin{array}{cc}
       \xymatrix@R=0.5em@C=1.8em{
    & A \cup \brC(A)   \ar@{-}[rd] \ar@{-}[ld] &  \\
    A \cup \o00 & & \brC(A) \\
     &\o00  \ar@{-}[ru] \ar@{-}[lu] &  \\
   }
    \end{array}
                                       \end{equation}
  is a diagram of submonoids, where  $(A \cup \o00) \cap \brC(A) = \o00$. Otherwise
  \begin{equation}\label{eq:2.5}
    \o00 \subset A \cup \o00 \subset \brC(A).
  \end{equation}
  All these submonoids are saturated in $V$.
\end{rem}

\begin{thm}\label{thm:2.9}
  Let $A$ be an archimedean class in $V$, and let $x$ be an element in $A$, whence $A =~ [x]_\aV$. Then
    \begin{equation}\label{eq:2.6}
      W = W(A) := W(x) :=  \{ y \in V \ds | \exists n \in \N : y \leq_V nx\}
    \end{equation}
    is the minimal SA-submonoid  of $V$ containing $A$, which is a saturated additive set (in short, a saturated SA-submonoid).
    It contains the  saturated submonoid $A \cup \brC(A)$. 
Both $W(A)$ and $A \cup \brC(A)$ are principal additive sets attracting the same entourage
 \begin{equation}\label{eq:2.6b}
 \brC(W(A)) = \brC(A \cup \brC(A)) = \brC(A).
\end{equation}
\end{thm}

\begin{proof} i) If $y_1 \leq_V n_1 x$ and  $y_2 \leq_V n_2 x$, then
$y_1 +y_2 \leq_V (n_1 + n_2) x$. Thus $W$ is a submonoid of~ $V$. If $y_1, y_2 \in V$ and $y_1 + y_2 \leq_V  m x$, then
$y_1 \leq_V m x$ and $y_2 \leq_V m x$ for some $m \in \N$. Thus   the submonoid $W$ is summand absorbing. Obviously it contains the set $A \cup \brC(A).$
\pSkip
ii)
Assume that $W'$ is an SA-submonoid of $V$ with $A = [x]_\aV \subset W'$. Let $y \in W$.
There exist $n \in \N$ and $z \in V$ for which  $y+z = nx \in W'$. Since $W'$ is summand absorbing, we infer that $y \in W'$. This proves that $W \subset W'$.
\pSkip
iii) From $A \subset A \cup \brC(A) \subset W(A)$ we conclude that $\brC(A) \subset \brC(A \cup \brC(A)) \subset \brC(W(A))$.
But  $W(A)$ is cofinal to $\N x$ by \eqref{eq:2.6}, now proven. Thus all three entourages are equal.
\end{proof}

Given an entourage $D \subset V$, we next  observe the flocks in the set  $T(D)$ of attractors of $D$,
 $$ T(D) = \{ A \in \Gm(V) \ds | \brC(A) = D\}, $$
via related saturated  submonoids of $V$. For each $A \in T(D)$ we have  already defined  the minimal SA-submonoid $W(A)$ of $V$ containing $A$, as described in Theorem \ref{thm:2.9}, and the saturated submonoid
$$ S(A):= A \cup \brC(A) = A \cup D$$
appearing in Diagrams \eqref{eq:2.4} and \eqref{eq:2.5}. Of course $S(A) \subset W(A)$.

Given $A \in \Gm(V)$, let $\Adw$ and $\Aup$ denote the sets of all $B,B' \in \Gm(V)$ with $B \leq_{\Gm(V)} A$ and
$A \leq_{\Gm(V)} B'$, respectively. Since $\Gm(V)$ is idempotent, we can read off from Lemma \ref{lem:1.6} that
$B \in \Adw$ \iff \ $B + _\Gm A =A$.
\begin{thm}\label{thm:2.10}
  $\Adw = \{ B \in \Gm(V) \ds| W(B) \subset W(A)\}. $
\end{thm}
\begin{proof} Given $A, B \in \Gm(V)$, we choose representatives $x$ and $y$ of $A$ and $B$, i.e.,
$A = [x]_\aV$ and  $B = [y]_\aV$. Then $B \in \Adw$ \iff \ $[y]_\aV + [x]_\aV = [x]_\aV$. This means that $[x+y]_\aV = [x]_\aV$, which   happens \iff \ there exists  $n \in \N$ with $y+x \leq_V nx$, and by Theorem \ref{thm:2.9} this means that $B \subset W(A)$. Thus,  $B \subset W(A)$ \iff \ $W(B) \subset W(A)$, since
 $W(B)$ is the minimal SA-submonoid of ~$V$ containing $B$.
\end{proof}

Given a fixed attractor $A$ of $D$, we look for  saturated submonoids larger than $S(A) = A \cup D$ that relate to the flock
$\Fl(A).$
\begin{prop}\label{prop:2.12}
The union $S_F(A)$ of all submonoids $S(A')$ with $A' \in \Fl(A)$ is
the minimal submonoid containing $D$ and all classes $A' \in \Fl(A)$.
\end{prop}

\begin{proof} Clearly $D \subset S_F(A)$. We are done, if we verify for any two classes $A_1', A_2' $ in $\Fl(A)$ that
 $[S(A_1')+ S(A_2')]_\sat \subset S_F(A)$.  We know that $B =  A_1' +_\Gm  A_2' \in \Fl(A)$,  and  furthermore that
$D + A_i' = A_i'$ for $i = 1,2$, and that $D+D =D$. It then follows that
\begin{align*}
S(A_1')+ S(A_2') & = ( D \cup A_1')+ (D \cup A_2') \\
& = D \cup ( D + A_1') \cup (D + A_2') \cup( A_1' + A_2')
\\ & =
D \cup A_1' \cup  A_2' \cup (A'_1+ A_2'),
\end{align*}
whence
\begin{equation*}\label{eq:2.8}
  \begin{array}{ll}
   [S(A_1')+ S(A_2')]_\sat  & = S(A_1') \cup  S(A_2') \cup (A_1' + A_2')_\sat \\[2mm]
   & =  S(A_1') \cup  S(A_2') \cup S(B).
  \end{array}
\end{equation*}
\vskip -5mm
\end{proof}

\begin{examp}\label{exmp:2.13}
Given a totally ordered set $J \subset \Gm(V)$, we give a rough account of the flocks contained in $J$.
Note that $J$ is a bipotent submonoid of $\Gm(V)$, i.e., $A+ B \in \{ A, B\}$ for all $A,B \in \Gm(V)$.  We choose a labelling
$J = \{A_\lm \ds | \lm \in \Lm  \}$, where $\Lm$ is a totally ordered index set and
$\lm < \mu$ \iff \ $A_\lm < _{\Gm(V)} A_\mu$. For each $\lm \in \Lm$ let $D_\lm = \brC(A_\lm)$. Then
\begin{equation*}\label{eq:2.9}
  \lm < \mu \dss \Rightarrow D_\lm \subset D_\mu, \text{ perhaps }  D_\lm = D_\mu.
\end{equation*}
The \textbf{maximal flock of $A_\lm$ in $J$} is
\begin{equation*}\label{eq:2.10}
  \Fl_J(A_\lm) := J \cap \Fl(A_\lm).
\end{equation*}
$\Fl_J(A_\lm)$ consists of  the attractors of $D_\lm$ in $J$ and all classes in $J \cap S_F(A_\lm)$. A further basic  fact is, that
\begin{equation*}\label{eq:2.11}
  \lm < \mu \dss \Leftrightarrow  W(A_\lm) \subsetneqq W (A_\mu)\dss \Leftrightarrow
  J \cap W(A_\lm) \subsetneqq J \cap  W (A_\mu).
\end{equation*}
If, say, $J$ has a maximal element $A_\om$ and a  minimal element $A_\al$, then while running down from ~$A_\om$ to $A_\al$ we arrive at a new flock precisely when the entourage $D_\lm$ changes, and  precisely when entering a smaller saturated submonoid $S_F(A_\lm) \supset D_\lm$ of $V$. 
Each maximal flock in $J$ is a convex subset of $J$.
\end{examp}
Next we  look for ``small'' flocks in an arbitrary additive monoid $V$. Given an archimedean class $A \in \Gm(V),$ it is evident by the definition of flocks (Definition \ref{def:2.6}.(b)) that, the intersection of any family of flocks  containing $A$ is again such a flock. Thus, there exists a  \textbf{unique minimal} flock containing $A$, denoted $\Flmin(A)$.

\begin{prop}\label{prop:2.14}
  Let $A \in \Gm(V)$ and let $D = \brC(A)$.
  \begin{enumerate} \ealph
    \item $\Flmin(A) = \Fl(A) \cap \Adw$.

    \item $\Flmin(A)$ is the downset of $A$ in $T(D)$.
  \end{enumerate}
\end{prop}
\begin{proof}
  (a): Let $\Phi = \Fl(A) \cap \Adw$. If $B \in \Phi$, then $\brC(B) = D$, since $B \in \Fl(A)$. Furthermore,  $B \leq_{\Gm(V)} A,$ and so $B + _\Gm A = A$, since $\Gm(V)$ is idempotent (cf. Lemma \ref{lem:1.6}). Thus,
  $\brC(A) = \brC(B) = \brC(A + B)$ by Theorem \ref{thm:1.3}. This proves that $A$ and $B$ are compatible.
  Hence, $\Phi$ is a flock. It is the minimal flock  containing $A$, since $\Phi \subset \Adw$.
  \pSkip
  (b): An immediate consequence of (a).
\end{proof}

We outline  some observations that  follow directly from  the definition of minimal flocks, without using Proposition \ref{prop:2.14}.

\begin{rem}\label{rem:2.15}
$ $
\begin{enumerate} \eroman
  \item If $A$ and $B$ are classes in the  same flock and  $A \leq_{\Gm(V)} B$, then $\Flmin(A) \subset \Flmin(B)$.
  \item If $A_1, A_2 \in \Fl(A)$, then $\Flmin(A_1 +_\Gm A_2)$ is the smallest flock containing $\Flmin(A_1)$ and
  $\Flmin(A_2)$.
\end{enumerate}
\end{rem}

Equipped with the concept of flocks,  we give a rough account of the essential and inessential archimedean classes in a principal additive subset $S$ of $V$, as defined above (Definition \ref{def:2.3}). Here we face a dichotomy  as follows.
Let $D= \brC(S)$. It is the entourage attracted by every essential class in $S$.

\begin{defn}\label{def:2.16}
If $A \in \Gm(V)$ is essential in $S$ and $B$
 is inessential in $S$, then, clearly, there are  two possibilities
 \begin{description}\dispace
   \item[Case I] $\brC(A+B) = \brC(A) = D$,
   \item[Case II] $\brC(A+B) = \brC(A) \supsetneqq D$.
 \end{description}
 We say in Case I, that $B$ is \textbf{controlled} in $S$, while in Case II, that $B$ is \textbf{excessive} in~ $S$.
We further say that $S$ is \textbf{controlled} (in $V$), if every inessential class of $S$ in controlled in ~$S$.
\end{defn}

This case distinction for $B$ does not depend on the choice of $A$ in $S$, cf. Definitions \ref{def:2.1}.(b) and \ref{def:2.6}.(b).

\emph{A memory aid.} The following view may be of help in understanding this terminology. Think of archimedean classes in $V $ as policemans. Every policeman $A$ is trained to supervise one district, namely  the entourage  $D = \brC(A)$. But he is not interested to report what happens in other districts. A flock $F$ is a family of policemans who work together in observing a given entourage $D$. There is a hierarchy in a flock $F$, given by the partial ordering $\leq_{\Gm(V)}$, restricted to $F$.  $F$ is a part of a unique maximal flock, that is $\Fl(A)$ for any $A \in F$. Thus every $A \in F$ can be supervised by one or several other officers in $\Fl(A)$, except the maximal members of $\Fl(A)$, if such members exist.

\begin{examp}\label{exmp:2.17}
  Let $A \in \Gm(V)$ and let $D = \brC(A)$. Assume that $A$ is not contained in $D$, whence disjoint from $D$. We introduce two principal additive subsets $S_1$ abd $ S_2$ of $V$:
  \begin{itemize}\dispace
    \item[--] $S_1:= \bigcup \Adw =$ the union of all classes in $\Adw$,
    \item[--] $S_2:= S_1 \cup D$.
  \end{itemize}
    Then $\brC(S_1) = \brC(S_2) = D$, i.e., $S_1$ and $S_2$ attract the same entourage $D$.
  They both have the same set of essential classes, that is  $\Adw \cap T(D)$. They also have the same set of controlled  inessential classes, namely  $\{ B \in \Adw \ds | B \subset D \}$. But $S_1$ is controlled,  while $S_2$ usually contains excessive inessential classes, the elements of $\Gm(D) \sm \Adw$.
  Here $\Gm(D)$ is viewed as a submonoid of $\Gm(V)$, cf. Proposition \ref{prop:2.7}.(b).
\end{examp}

\begin{examp}\label{exmp:2.18}
 Let $Z \in \Gm(V)$, let $D = \brC(Z)$, and assume that $Z \subset D$.
 This is Case I in Dichotomy \ref{dich:2.2} above (where we now write $Z$ instead of $A$). Then
 $$ T(D) = \Zup = \{ B \in \Gm(V) \ds| Z+_\Gm B = B\} $$
 is a flock. (Note that $\Zup \subset D$, since $D$ is summand absorbing in $V$, cf. Proposition \ref{prop:2.7}.(b).)
 Thus, $T(D)$ is the \textbf{unique maximal flock} $\Fl(Z)$ with member $Z$, and $\bigcup \Zup$ is the unique maximal principal additive set containing the class $Z$.

We say that the entourage $D$ is \textbf{centered} and call $Z$ the \textbf{center} of $D$ (as above), and also name $\Fl(Z)$ \textbf{centered}. $\Zup \sm \{ Z\}$ is the set of inessential classes of $\bigcup \Zup = \Fl(Z)$.
\end{examp}
We return to principal additive sets in general.
\begin{defn}\label{def:2.19} Given a principal additive subset $S$ of $V$, and $D = \brC(S)$, we define the ~set
\begin{equation*}\label{eq:2.12}
  F_S := \{ A \in \Gm(V) \ds | A \subset S, \ \brC(S) = D\}
\end{equation*}
of essential classes in $S$.
Conversely, for a flock $F$ in $V$ attracting $D$ we define  the set
\begin{equation*}\label{eq:2.13}
S_F:= \bigcup F = \text{the union of all } A \in F.
\end{equation*}
We call $F_S$ the \textbf{flock associated to} $S$, and say that $S_F$ the \textbf{principal additive set associated to $F$}.
\end{defn}
It is evident from Definitions \ref{def:2.3}.(b) and \ref{def:2.6}.(b), that in this  way the flocks in $V$ correspond uniquely to the principal additive subsets of $V$.

\begin{rem}\label{rem:2.20} Given a principal additive set $S$ of $V$ and an essential class $A$ in $S$, it is obvious that for any $B \in \Gm(V)$ with $B <_{\Gm(V)} A$ and $\brC(A) = \brC(B) = \brC(S)$ the set $S' = S \cup B$ is again
principal additive and $\brC(S') = \brC(S)$. Thus we may add, and of course also omit, classes ``at the bottom'' of $S$ without changing the entourage $\brC(S)$ and the sets of controlled and excessive inessential classes in $\brC(S)$.
\end{rem}

\begin{defn}\label{def:2.21}
We say that a principal additive subset $S$ of $V$ is \textbf{grounded}, if
$$ F_S = \Fdw \cap T(D), \qquad D:= \brC(S).$$
This means that $S$ cannot be enlarged by adding any class $B$ as above.
\end{defn}

For example, both sets $S_1$ and $S_2$ appearing in Example \ref{exmp:2.17} are grounded. Clearly, the set $S_F$
associated to a maximal flock $F$ is grounded.

\section{The minimal archimedean classes}\label{sec:3}
Within archimedean classes a particular interest is devoted to minimal classes.

\begin{defn} \label{def:3.1} $ $
 \begin{enumerate}\ealph
   \item We call an archimedean class in $V$ \textbf{minimal}, if there is no class $B \in \Gm(V)$ with
   $\brC(A) = \brC(B)$ and $B < A$ in $\Gm(V)$.
   This means that $\Flmin(A) = \{ A\}$, cf. Proposition~ \ref{prop:2.14}.
   \item Starting with an entourage $D$, we denote by $T(D)_{\min}$ the set of all minimal classes ~$A$ with $\brC(A) = D$.
 \end{enumerate}
 \end{defn}
 Of course, it may well happen that  $T(D)_{\min}$ is empty. Assume  that $T(D)_{\min} \neq \emptyset$, and let  $A \in T(D)_{\min}$. We explore the minimal summand absorbing submonoid $W(A)$ of $V$ containing ~ $A$. By Dichotomy \ref{dich:2.2},  there are two cases: either $A \subset \brC(A) = D$ or $A$ is disjoint from $D$.

 In the first case  $T(D)_{\min}$ consists of a unique class $A = Z$, and $W(Z) = D$, i.e., $D$
is centered, cf. Remark \ref{rem:2.8}. Then $D$ admits a unique maximal flock of attracting classes, namely $\Zup = T(D)$.

We now assume that $D$ is an entourage in $V$ \emph{which is not centered}. We pick a class $A \in  T(D)_{\min}$ and choose an indexing
\begin{equation*}\label{eq:3.1}
  \Fl(A) \cap  T(D)_{\min} = \{ A_i \ds| i \in I\}
\end{equation*}
of all minimal archimedean classes which are relatives of $A$. Thus, all $A_i$ are disjoint from~ $D$
(and, of course, also $A_i \cap A_j = \emptyset$ for all $i \neq j$).

\begin{thm}\label{thm:3.2}
If $A$ is minimal and $\brC(A) = D$ is not centered, then $W(A) $ is the disjoint union of the minimal classes in
$\Fl(A)$ and $D$, i.e.,
$$ W(A) = \dot {\bigcup_{i \in I}}\; A_i \; \dot \cup \; D.$$

\end{thm}
\begin{proof}
  We conclude from $A_i \subset W(A)$ that $W(A_i) \subset W(A)$. But $A_i \not \subset D$, and thus $D \subsetneqq W(A_i)$.
  Since $A$ is minimal, it follows that $W(A_i) = W(A)$ for all $i \in I$. Given an essential class ~$B$ of $W(A)$, we see by  the same way that $B$ is minimal in $\Fl(A)$, and hence $B= A_i$ for some $i \in I.$
\end{proof}

The case $|I| =1$ deserves special interest.
\begin{defn}\label{def:3.3}
We call an archimedean class $A$ in $V$ \textbf{tight}, if $A \cup \brC(A) = W(A)$.
\end{defn}

\begin{cor}\label{cor:3.4}
  An archimedean class $A$ in $V$ is tight \iff \  either $A$ is central, i.e., $A \subset \brC(A)$, or $A$ is not centered and
  $\Fl(A) = \{ A\}$.
\end{cor}
\begin{proof}
  Clear by Theorem \ref{thm:3.2} and the discussion on central classes preceding it.
\end{proof}

\begin{schol}\label{schol:3.5}
  We  explicitly observe the condition that  an archimedean class in $V$ is tight, leaving aside the case that $\brC(A)$ is centered. Given an element $x \in V$ with $A = [x]_\aV$, then  $x \notin \brC_\om(x) = \bigcup_{n \in \N} \brC(nx)$ (Case II in Dichotomy \ref{dich:2.2}). Thus, $A$ is tight \iff \ (cf. ~ \eqref{eq:2.6})
  $$ W(A) = W(x) = \{y \in V \ds | \exists n \in \N : y \leq_V nx \} $$
  is the disjoint union of $[x]_\aV$ and $\brC_\om(x)$. This means that for any $y \in V$ the following holds.
  \begin{equation*}\label{eq:3.2}
\begin{array}{llrl}
  \exists n \in \N : y \leq_V nx & \ds\Rightarrow & \text{ either } & \exists m \in \N: x \leq_V m y \\
  &&   \text{ or } & \exists n' \in \N: y+n'x = n'x.
  \end{array}
   \end{equation*}
  \end{schol}

Starting from Theorem \ref{thm:3.2}, we aim  to establish a natural hierarchy in a given non centered maximal flock. To this end, as a preparation,  we state a kind of duality for certain submonoids of $V$. (It could have been proved much earlier.)

\begin{prop}\label{prop:3.6}
  Assume that $W$ is a submonoid of $V$, and let $U= (V \sm W) \cup \o00$.
  \begin{enumerate}\ealph
    \item $W$ is SA in $V$ \iff \ $U$ is a submonoid of $V$ \iff \  $V \sm W$ is closed under addition.
    \item $W$ is saturated in $V$, i.e., $W$ is a union of classes of the equivalence relation $\equiv_V$, \iff \ this holds for $U$.
    \item If $W$ is SA in V, then the quasiordering $\leq_V$ restricts to the quasiorderings $\leq_U$ on $U$ and $\leq_W$ on $W$.
   \item The ordering $\leq_\aV$ on $\Gm(V)$ restricts to $\leq_\aU$ and $\leq_\aW$ on $\Gm(U)$ and $\Gm(W)$.
  \end{enumerate}
\end{prop}
\begin{proof}
  Claims (a) and (c) are immediate consequences of the SA-condition for submonoids of ~$V$, and (b) is evident. Claim (d) follows from claim (c), since that relation $\leq_\aV$ on $V$ is a derivate of $\leq_V$:
  $$ y \leq_\aV x \dss\Leftrightarrow \exists n \in \N: y \leq_V nx,$$
  and the same holds for $W$ and $U$.
\end{proof}

We will also need a general transitivity lemma for flocks.

\begin{lem}\label{lem:3.7}
  Assume that $F_1$ is a flock in an additive monoid $V$ attracting the submonoid~ $D$, and that $F_2$ is a flock in $(V \sm D) \cup \o00$ attracting $S_{F_1}$, cf. Definition \ref{eq:2.13}. Then $F_1 \cup F_2$ is a flock in~ $V$ attracting $D$.
\end{lem}
\begin{proof}
  Given classes $A \in F_1$ and $B \in F_2$, we choose elements $x,y \in V$ with $A = [x]_\aV$ and $B = [y]_\aV$. Then $x+y \equiv_\aV y $ and $$A+_{\Gm(V)} B = [x+y]_\aV = [y]_\aV =B.$$
  This proves that $A$ and $B$ are compatible (= relatives) in $V$ (Definition \ref{def:2.6}). Of course all members of $F_1$ and all members of $F_2$ are compatible in $V$, and so all members of $F_1 \cup F_2$ are compatible in $V$. Thus, $F_1 \cup F_2$ is a flock. It attracts $D$, since $F_1$ attracts $D$.
\end{proof}
\begin{thm}\label{thm:3.7}
  Assume that $F$ is a maximal flock in $V$ with $\brC(F) = D$ disjoint from every $A \in F$ (Case II in Dichotomy \ref{dich:2.2}), and that the set $\Fmin$ of minimal elements of $F$ in $\Gm(V)$ is not empty, with the consequence that $\Adw$ contains a minimal element for any $A \in F$. Let ~$W(F)$ denote the minimal SA-submodule of $V$ containing $F$.
  \begin{enumerate}\ealph
  \item
   Then
   $$V_1(F) := [W(F) \sm D] \cup \o00 = W(F) \sm (D \sm \o00)$$ is a saturated additive submonoid of ~$V$.

  \item The quasiordering $\leq_V$ on $V$ restricts to the quasiordering $\leq_{V_1}$ on $V_1(F)$.

  \item The ordering $\leq_\aV$ on $\Gm(V)$ restricts to the ordering $\leq_{\aV_1}$ on $\Gm(V_1(F))$. \\
  \{N.B. $\Gm(V_1(F)) \subset \Gm(V)$ by Proposition \ref{prop:2.7}.\}

   \item If $F \neq \Fmin$, then $$D_1(F) := \big(\bigcup \Fmin \big) \cup \o00$$ is an entourage in $V_1(F)$, and $ F \sm \Fmin$ is a maximal flock in $V_1(F)$ attracting $D_1(F)$.
\end{enumerate}
\end{thm}

\begin{proof}
  We may replace  $V$ by $W(F)$, since $W(F)$ is SA in $V$, and thus  $\leq_V$ restricts to ~ $\leq_{W(F)}$, and furthermore all claims hold for $W(F)$. Now assertions (a)--(c) are clear by Proposition~ \ref{prop:3.6}. Part  (d) is covered by
  Theorem \ref{thm:3.2}.\end{proof}
Next we  iterate the construction of $V_1(F)$ and $D_1(F)$ in Theorem \ref{thm:3.7}, assuming again that~ $F$ is a maximal flock in $V$, where  $\Fmin \neq \emptyset$ and $\brC(F) = D$ are disjoint from every~ $A \in F$.

\begin{construction}\label{cons:3.8}
  We define the sequence of subsets $F_1, F_2, \dots $ of $F$,
  \begin{equation}\label{eq:3.3}
    F_1 := \Fmin, \quad F_2 = (F \sm F_1)_{\min}, \quad F_3 = (F \sm (F_1\cup F_2))_{\min}, \quad  \dots
  \end{equation}
  which  gives three possible decompositions:
  \begin{description}
    \item[Case I]
  $$ F = F_1 \ds{\dot \cup} \cdots \ds{\dot \cup} F_d \ ,$$
  \item[Case II]
  \begin{equation*}\label{eq:3.4}
  F = \dot {\bigcup_{i \in \N}}\; F_i \ ,
  \end{equation*}
\item[Case III]
\begin{equation*}\label{eq:3.5}
  F = \dot {\bigcup_{i \in \N}}\; F_i \ds{ \dot \cup} G = H \ds{ \dot \cup} G\ ,
  \end{equation*}
where $H = \dot \bigcup_{i \in \N} F_i $ and  $G = F \sm \dot \bigcup_{i \in \N} F_i \neq \emptyset$.
\end{description}
In consequence of Theorem \ref{thm:3.7} we obtain  the submonoids $V_1(F), V_2(F), \dots$
and $D_1(F)$, $D_2(F), \dots$ of $W(F)$ as follows
$$ V_1(F) = (W(F)\sm D) \cup \o00.$$
If $F \neq F_1$, then $F \sm F_1$ is a maximal flock in $V_1(F)$, attracting $D_1(F)$.
If $F \neq F_1 \cup F_2$, then $F \sm (F_1 \cup F_2)$ is a maximal flock in $$V_2(F) = [V_1(F) \sm D_1(F)] \cup \o00,$$ attracting $D_2(F)$, etc .
\end{construction}

\begin{thm}\label{thm:3.9} $ $
\begin{enumerate} \ealph
  \item Let $F = F_1 \ \dot \cup \cdots \dot  \cup \ F_d$ (Case I). Then $F$ contains a \textbf{unique tight} archimedean class $\tau = \tau_F$. It is the unique maximal element of $F$ as a subset of $\Gm(V)$, $F_d = \{\tau \}$, and $F = \tau^\downarrow \cap T(D)$.
  \item Let $F = \bigcup_{i \in \N} F_i $ (Case II), and let $S_i = S_{F_i}$ in $V_i(F)$, cf. Definition \ref{eq:2.13}. Then $\brC(S_i) = D_i(F)$, $S_F = \bigcup_ {i \in \N} S_i$ (in $W(F)$ or in $V$), and $W(F) = S_F \ \dot \cup \ D$.
  \item Let $F = H \ds{\dot\cup} G$, where  $ H= \bigcup_{i \in \N} F_i $ and   $G \neq \emptyset$ (Case III). Then
      $G = F \sm  H$ is a \textbf{maximal} flock  in $W(F) \sm S_H$. 
\end{enumerate}

\end{thm}

\begin{proof}
 \pSkip
  (a):  The flock $F_d$ cannot contain two different classes $A_1$ and $A_2$, since this would imply that  there is a class
  $A_1 +_{\Gm} A_2$ in $F$  which is not contained in $F \sm (F_1 + \cdots + F_d)$. Thus, ~ $F_d$ consists of a dingle  class $\tau$. It is tight by Corollary \ref{cor:3.4}, due to the definition of the ~$F_i$, whence $F = \tau^\downarrow \cap T(D)$.

   \pSkip
  (b): Obvious by the one-to-one correspondence of flocks and principal additive sets, stated following Definition \ref{def:2.19}.

 \pSkip
  (c): Again obvious by this one-to-one correspondence, except for the claim that the flock~ $G$ is maximal in $W(F) \sm S_H$. But, if there would exist a flock $G' \supsetneqq G$ in $W(F) \sm S_H$, then by Lemma~ \ref{lem:3.7} the flock
  $\bigcup_{i \in \N} F_i \ds\cup G'$ would be larger   than the maximal flock $F$ -- a contradiction. Thus,~ $G$ is maximal in $W(F) \sm S_H$.
\end{proof}
In the following we view the maximal flock $F$ as a subset of the idempotent monoid $\Gm(V)$. In Cases I and II of Construction \ref{cons:3.8},
we define a \textbf{height function} $h: F \to \N$ by the rule
\begin{equation}\label{eq:3.6}
  h(A) := k \dss{\Leftrightarrow} A \in F_k.
\end{equation}
In Case I the range of this function is $\{1, \dots, d \}$, while in Case II it is $\N$. In Case III we define a height function $h: H \to \N$ by the same rule \eqref{eq:3.6}.
\begin{thm}\label{thm:3.11} $ $ Assume that $A$ and $B$ are archimedean classes in $F$ or in $H$ respectively.
Then 
\begin{enumerate} \ealph
\item $\max(h(A),h(B)) \leq h(A+_\Gm B)$.

\item $\max(h(A),h(B)) < h(A+_\Gm B)$ \iff \  $A$ and $B$ are incomparable in the idempotent monoid ~$\Gm(V)$.

\item  $ h(A+_\Gm B) \leq h(A) + h(B)$.

\item $A < B \Rightarrow h(A) < h(B)$.

\end{enumerate}
\end{thm}
\begin{proof}
   By the general Lemma \ref{lem:1.6}.(a), we have $A < B$ \iff \ $A +_\Gm B = B$, and $B < A$ \iff \
  $A +_\Gm B = A$.   Otherwise
  $A +_\Gm B >  A$ and  $A +_\Gm B > B$, again by this lemma. This covers claims  (a) and (b).

\pSkip (c): The claim can be settled by induction on $h(B)$. If $B$ and $C$ are classes in $F$, respectively in $H$, then  $h(B+_\Gm C) = h(B)$, if $C < B$,  and  $h(B+_\Gm C) = h(B) +1$ otherwise. From this we infer that  $h(A+_\Gm B) \leq h(A) + h(B)$, since $B$ is a sum of elements of height $1$ in the submonoid~ $F$,  respectively $H$, of $\Gm(V)$.

\pSkip
 (d): Let $A < B$, where  $A= F_k$ and  $B= F_\ell$. If $k = \ell$, the classes $A$ and $ B$ would be incomparable, since then both are minimal in the flock $F \sm \bigcup_{i <k } F_i$. Thus $k \neq \ell$. We conclude by \eqref{eq:3.6} that $k = h(A) < \ell = h(B)$.
\end{proof}

\begin{defn}\label{def:3.12} In Case III of Construction \ref{cons:3.8} we call $H$ the \textbf{initial part} and call  $G = F \sm H$ the \textbf{transfinite  part} of $F$. We denote these sets more elaborately by
$H_F$ and ~$G_F$, if necessary.
In Case I and II we formally put $H_F = F$.
\end{defn}
It may well happen that $G_{\min} =  \emptyset,$ but otherwise we can repeat the procedure  with
$(G, G_H)$ in $W(G) = W(F)$ instead of $(F,D)$, to  obtain a second $\N$-valued  height function $h_G$ on $G$ or on the initial part $H_G$ of $G$.

We extend the height function $h: H_F \to \N$ to a height function
$$ h': H_F \ds{\dot \cup} H_G \To \On \sm \00$$
taking  values in the subset $\N \cup (\om + \N)$ of non-limit ordinals $\al $ with $1 \leq \al \leq \om + \om$, defined by the rule
\begin{equation*}\label{eq:3.7}
h'(A) = \left\{   \begin{array}{ll}
    h(A), & \text{ for } A \in H_F, \\[1mm]
    \om + h_G(A), & \text{ for } A \in H_G.
    \end{array} \right.
\end{equation*}
\begin{thm}\label{thm:3.12} For any $A,B\in H_F \ \dot \cup \ H_G$
\begin{enumerate} \ealph
\item $\max(h'(A),h'(B)) \leq h'(A+_\Gm B)$.

\item $\max(h'(A),h'(B)) < h'(A+_\Gm B)$ \iff \ $A$ and $B$ are incomparable in the idempotent monoid $\Gm(V)$.

\item  $ h'(A+_\Gm B) \leq h'(A) + h'(B)$, if $A\in H_F$ or $B  \in H_F$.

\item $A < B \Rightarrow h'(A) < h'(B)$.

\end{enumerate}
\end{thm}
\begin{proof}
  (a) and (b) hold essentially by the same arguments as in the proof of Theorem \ref{thm:3.11}.
  \{$A < B \Leftrightarrow A + _\Gm B = B$, etc.\}

\pSkip
  (c): Holds by the induction argument used for Theorem \ref{thm:3.11}.(c), since here $A$ or $B$ is a sum of classes in $F_{\min}$.

\pSkip
  (d): If both $A$ and $B$ are in $H_F$ or in $H_G$, then $h'(A) < h'(B)$ by Theorem \ref{thm:3.11}.(a). Otherwise $A< B$ forces that $A \in H_F$, $B \in H_G$, $h'(A) < \om$, and  $h'(B) > \om$. Thus $h'(A) < h'(B)$.
\end{proof}

\end{document}